\documentclass{amsart}  
\usepackage[T1]{fontenc}
\usepackage{listings}
\usepackage{mathrsfs}
\usepackage{graphicx}
\usepackage{bbm}
\usepackage{comment}
\usepackage{amssymb}
\usepackage{enumerate}
\newtheorem{theorem}{Theorem}[section]
\newtheorem{lemma}[theorem]{Lemma}
\newtheorem{cor}[theorem]{Corollary}
\theoremstyle{definition}

\newtheorem{rem}[theorem]{Remark}

\numberwithin{equation}{section}

\newcommand\F{\mathcal{F}}
\newcommand\E{\mathbb{E}}
\newcommand\R{\mathbb{R}}

\begin{document}

\title{On the diameter of the stopped spider process}

\author{Ewelina Bednarz}
\address{Faculty of Mathematics, Informatics and Mechanics, University of Warsaw, Banacha 2, 02-097 Warsaw, Poland}
\email{ewela.bednarz@gmail.com}

\author{Philip A. Ernst}\thanks{}
\address{Department of Statistics, Rice University, Houston, TX 77005, USA.}
\email{philip.ernst@rice.edu}

\author{Adam Os{{\k e}}kowski}\thanks{}
\address{Faculty of Mathematics, Informatics and Mechanics, University of Warsaw, Banacha 2, 02-097 Warsaw, Poland}
\email{A.Osekowski@mimuw.edu.pl}



\begin{abstract}
We consider the Brownian ``spider process'', also known as Walsh Brownian motion, first introduced in the epilogue of Walsh (\cite{Walsh1978}). 
The paper provides the best constant $C_n$ for the inequality
$$ \E D_\tau\leq C_n \sqrt{\E \tau},$$
where $\tau$ is the class of all adapted and integrable stopping times and $D$ denotes the diameter of the spider process measured in terms of the British rail metric. The proof relies on the explicit identification of the value function for the associated optimal stopping problem.
\end{abstract}

\maketitle

\section{Introduction.} 
We consider the Brownian ``spider process'', also known as Walsh Brownian motion, as first introduced in the epilogue of Walsh 1978 (\cite{Walsh1978}). Early constructions of Walsh's Brownian motion were given by L.C.G. Rogers (\cite{Rogers}) using resolvents, by Baxter and Chacon (\cite{Baxter}) using infinitesimal generators, and by T.S. Salisbury (\cite{Salisbury}) using excursion theory. The 1989 work of Barlow, Pitman, and Yor (\cite{BPY}) considered the construction of Walsh Brownian motion as a process living on $n\ge 1$ rays meeting at a common point (reminiscent of a spider), and it is this construction that we shall exploit in the present paper.\\
\indent The purpose of this paper is to solve an optimal stopping problem for the spider process which we shall give in equation \eqref{optimal_stopping} below. Before revealing the problem, we begin with some background and necessary definitions. Conceptually, the construction of the spider process is motivated by the fundamental observation that standard one-dimensional Brownian motion can be viewed as an absolute value of itself, each of whose excursions is assigned a random sign.  Following the construction in \cite{BPY,DS}, the spider process, with $n \geq 3$ rays emanating from the origin, may be viewed as the extension of the above observation to an $n$-valued sign. More precisely, for a given positive integer $n$, consider the collection of $n$ rays $ R_k=\{e^{2\pi i k/n}t\,:\,t\geq 0\}$, $k=1,\,2,\,\ldots,\,n,$ on the plane.   Let $\theta=(\theta_m)_{m\geq 0}$ be the sequence of independent ``complex signs'', i.e., the family of random variables with the uniform distribution on $\{e^{2\pi i k/n}:k=1,2,\ldots,n\}$. 
Assume further that $B=(B_t)_{t\geq 0}$ is a Brownian motion independent of $\theta$ and let $e=(e_t)_{t\geq 0}$ be the associated excursion process (see Chapter XII in \cite{RY}). The set of excursions is countable and hence it can be ordered by the set of natural numbers. The spider process $S$ is then given by $ S_t=\theta_{m(t)}|B_t|$, where $m(t)$ is the number of the excursion of $B$ which straddles $t$ (see \cite{RY}, p. 488). From this definition, we see that for the case $n=1$, the spider process reduces to reflecting Brownian motion $|B|$; the case $n=2$ corresponds to standard Brownian motion.\\ 
\begin{figure}[htbp]
\begin{center}
\includegraphics[scale=0.3]{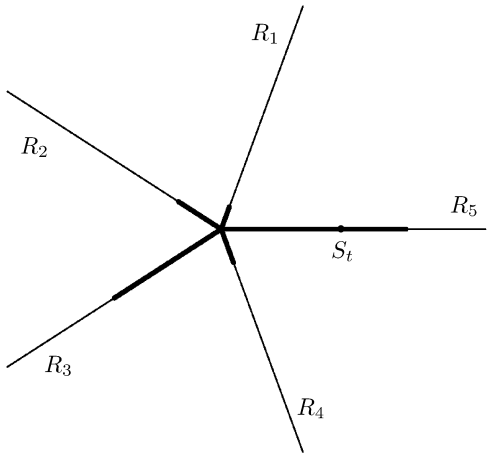}
\caption{A depiction of the spider process with $n=5$. The bold segments indicate the points already visited by $S$ up to time $t$. The longest rib lies on the ray  $R_5$. The the second longest rib lies on $R_3$. }
\end{center}
\end{figure}
\indent The optimal stopping problem in \eqref{optimal_stopping} is motivated by the development of optimal bounds for the expected ``size'' (as defined in equations \eqref{size1} and \eqref{size2} below) of the stopped spider process. For a given $t\geq 0$ and $\omega\in \Omega$, let $T_t(\omega)$ denote the trajectory up to time $t$: 
\begin{equation}\label{trajectory}
T_t(\omega)=\{S_s(\omega):0\leq s\leq t\}.
\end{equation}
Moreover, let $d_t(\omega)$ stand for the sum of the deviations of $T_t(\omega)$ along the rays, i.e.,
\begin{equation}\label{diameter}
 d_t(\omega)=\sum_{k=1}^n |T_t(\omega)\cap R_k|.
\end{equation}
We shall sometimes refer to these deviations as the ``ribs'' of $S$.\\
\indent The late Lester Dubins considered the problem of designing a stopping time to maximize the coverage
of Brownian motion on the spider for a given expected time (see \cite{E}). That is, he sought to find the best constant $c_n$ such that
\begin{equation}\label{Dubins_question}
 \E d_\tau\leq c_n\sqrt{\E \tau}
\end{equation}
for any integrable stopping time $\tau$ of $S$  (i.e., measurable with respect to the filtration generated by the spider process). To the best of our knowledge, this question has been answered only in the cases $n=1$ and $n=2$. When $n=1$, $S=|B|$ and hence \eqref{Dubins_question} becomes
$$ \E \sup_{0\leq t\leq \tau}|B_t|\leq c_1 \sqrt{\E \tau}.$$
Dubins and Schwarz \cite{DS} proved that the value $c_1=\sqrt{2}$ is optimal; one may also consult Dubins, Gilat and Meilijson \cite{DGM} for an alternative approach.
 For $n=2$, the spider process reduces to the Brownian motion and we define $d$ to be the difference of the running maximum and the running minimum, namely $$d_t=\sup_{0\leq s\leq t}B_s-\inf_{0\leq s\leq t}B_s.$$ In \cite{DGM}, the authors proved that the optimal choice for $c_2$ is equal to $\sqrt{3}$. For $n\geq 3$, a tempting conjecture is that $c_n=\sqrt{n+1}$, but this appears to not be so, at least for $n=3$ (cf. \cite{E}). 

 In this paper, we shall study a version of Dubins' question in which the coverage or size of the spider process is measured differently. In this formulation, we shall replace $d_t$ by the diameter $D_t$ with respect to the British rail metric: that is, for $n=1$, we have
\begin{equation} \label{size1}
D_t(\omega)=|T_t(\omega)\cap R_1|.
\end{equation}
For $n\geq 2$, we have
\begin{equation} \label{size2}
D_t(\omega)=\max\Big\{\big(|T_t(\omega)\cap R_j|+|T_t(\omega)\cap R_k|\big)\,:\,j,\,k\in \{1,2,\ldots,n\}, j\neq k \Big\}.
\end{equation}
 In other words, $D_t(\omega)$ is given as in \eqref{diameter}, but only one or two largest summands are taken into account (depending, respectively, on whether $n=1$ or $n=2$). In the simpler case that $n\in \{1,2\}$, then  $d_t=D_t$, so the optimal constant $C_n$ in the inequality
\begin{equation}\label{main_inequality}
\E D_\tau\leq C_n \sqrt{\E \tau},
\end{equation}
equals $\sqrt{2}$ for $n=1$ and $\sqrt{3}$ for $n=2$. The key purpose of the present paper is to identify the optimal value of $C_n$ for $n\geq 3$. We preview the paper's main result as Theorem \ref{main_theorem} below, which shall be proved in Section \ref{sec4}.

\begin{theorem}\label{main_theorem}
For $n\geq 3$, the best constant in \eqref{main_inequality} is given by
$$ C_n=2\sqrt{U(0,0,0)},$$
where $U(0,0,0)$ is defined in Corollary \ref{const2} below.
\end{theorem}

A few remarks about the general strategy for proving our main result are now in order. By a straightforward time-homogeneity argument, it suffices to find the optimal constant $\kappa_n$ in the inequality $$\E D_\tau- \E \tau\leq \kappa_n. $$ This directly leads to the optimal stopping problem
\begin{equation}\label{optimal_stopping}
 \mathbb{U}=\sup \E(D_\tau-\tau),
\end{equation}
where the supremum is taken over all integrable stopping times $\tau$ of $S$. 
Given an optimal stopping problem, one generalizes the problem to the case in which the underlying Markov process is allowed to start from an arbitrary point from the associated state space. As a result, the corresponding value $\mathbb{U}$ extends to the value (or ``reward'') function on the entire state space. This object has many structural properties which in many cases enables its explicit identification (which, in turn, yields the solution to the initial optimal stopping problem). To find the reward function, one typically exploits one of the following two strategies:

\smallskip

(A) Using Markovian arguments, we may write a system of differential equations that the value function must/should satisfy. Then, applying analytic arguments and exploiting homogeneity (if there is any) in the problem structure, we  attempt to solve such a system and come up with the ``right'' function.

\smallskip

(B) We attempt to guess the the optimal strategy. To do so, we compute the value function from the very definition (by specifying, for each starting point, the corresponding optimal stopping time).

\smallskip

Sometimes, a successful solution requires a clever combination of both (A) and (B). It should be emphasized that typically  both these approaches yield only the \emph{candidate} for the reward function (during the search and the construction for the reward function, one usually exploits a number of guesses and/or some additional assumptions, which are not a priori guaranteed). Next, having completed this informal step (i.e., having found the candidate), one proceeds to rigorous proof and checks the excessiveness and optimality of the constructed function. If both excessiveness and optimality hold, then the candidate coincides with the value function and the optimal stopping problem is solved.

In solving the optimal stopping problem in \eqref{optimal_stopping}, we shall utilize both of strategies (A) and (B). We will also need a number of novel arguments; in particular, in order to reduce dimensionality and represent the spider process in terms of a relatively simple Markovian structure, we shall employ a skew Brownian motion with jumps. Furthermore, by a certain translation property and an appropriate reduction trick (Section \ref{sec2}), we shall see that the analysis of optimal stopping problem in \eqref{optimal_stopping} will heavily depend on the solution of a related auxiliary two-dimensional optimal stopping problem \eqref{auxiliary_stopping_prob} for a standard Brownian motion.

The remainder of this paper is organized as follows. Section \ref{sec2} is concerned with the analysis of the aforementioned auxiliary stopping problem \eqref{auxiliary_stopping_prob}. For the sake of completeness, we also present the solution to \eqref{optimal_stopping} for the cases $n=1$ and $n=2$. Although  the solution for both these cases has already appeared in the literature, we shall find that their analyses give helpful intuition about the optimal strategy in the more general case. Section 3 is devoted to the construction of the candidate for the value function for $n\geq 3$. It is the most technically innovative part of the paper; our efforts shall include the aforementioned reduction as well as a combination of arguments from methods (A) and (B) above. Section \ref{sec4} concludes the manuscript by first proving that the constructed candidate coincides with the desired value function, and then showing how the optimal stopping problem in \eqref{optimal_stopping} leads to the proof of our main result (Theorem \ref{main_theorem}).

\section{Preparation.}\label{sec2}

\subsection{A related optimal stopping problem} We begin with a problem, which itself is not new (see, for example, \cite{E}), but whose analysis will be quite helpful in for the present paper. For the sake of convenience and clarity, the presentation is divided into a few intermediate parts.

\smallskip

\emph{Step 1.} Suppose that $X=(X_t)_{t\geq 0}$ is a standard, one-dimensional Brownian motion and denote by $Y=(Y_t)_{t\geq 0}$ the associated one-sided maximal function, i.e., $Y_t=\sup_{0\leq s\leq t}X_s$ for $t\geq 0$. 
Consider the optimal stopping problem
\begin{equation}\label{auxiliary_stopping_prob}
 \mathbb{V}=\sup \E \left(Y_\tau-\tau\right),
\end{equation}
where the supremum is taken over all integrable stopping times $\tau$ of $X$. By Wald's identity, the time variable can be removed: the above supremum equals $
 \mathbb{V}=\sup \E \left(Y_\tau-X_\tau^2\right)$, thus giving rise to the optimal stopping problem for the Markov process $(X,Y)$. We distinguish the associated state space $$\mathcal D=\{(x,y):y\geq x\vee 0\}$$ and introduce the gain function $G:\mathcal D\to \R$ given by  $G(x,y)=y-x^2$. We then have the identity
\begin{equation}\label{stop_n=1}
 \mathbb{V}=\sup \E G\left(X_\tau,Y_\tau\right),
\end{equation}
which fits into the general framework of the theory of optimal stopping (cf. \cite{PS}). As mentioned in the previous section, a successful treatment of \eqref{stop_n=1} requires the generalization of the problem to the case in which the process $(X,Y)$ starts from an arbitrary point in the state space $\mathcal {D}$. This is standard: one first extends the process $(X,Y)$ to a Markov family on $\mathcal D$, introducing the  family of initial distributions  $(\mathbb{P}_{x,y})_{(x,y)\in \mathcal D}$, given by the requirement that $$\mathbb{P}_{x,y}(X_0=x,Y_0=y)=1$$ for all $(x,y)\in \mathcal{D}$. Next, one defines the associated value function
\begin{equation}\label{stop_n=1'}
 \mathbb{V}(x,y)=\sup \E_{x,y}G(X_\tau,Y_\tau),\qquad (x,y)\in \mathcal D,
\end{equation}
where the supremum is taken over all $\mathbb{P}_{x,y}$-integrable stopping times $\tau$ of $X$. Alternatively, one can define $\mathbb{V}(x,y)$ using a single probability measure $\mathbb{P}_{0,0}$ by
\begin{equation}\label{0,0}
 \mathbb{V}(x,y)=\sup\E_{0,0}G\left(x+X_\tau,\left(x+\sup_{0\leq s\leq \tau}X_s\right)\vee y\right),
\end{equation}
where the supremum is taken over all $\mathbb{P}_{0,0}$-integrable stopping times $\tau$ of $X$. 

\smallskip

\emph{Step 2.} Following the usual approach from general optimal stopping theory (cf. \cite{PS}), we split the state space $\mathcal {D}$ into two sets, the continuation region $C$ and the instantaneous stopping region $D$. They are given, respectively, by
$$ C=\big\{(x,y)\in \mathcal{D}\,:\, \mathbb{V}(x,y)>G(x,y)\big\}, \quad  D=\big\{(x,y)\in \mathcal{D}\,:\, \mathbb{V}(x,y)=G(x,y)\big\}.$$
Thus, to solve \eqref{stop_n=1'} (and hence also \eqref{stop_n=1}), one needs to identify the shape of the continuation region and the formula for $\mathbb{V}$ on this set. Having done that, the optimal stopping time is given by
\begin{equation}\label{opt_stop}
 \tau=\inf\{t\geq 0:(X_t,Y_t)\in D\}.
\end{equation}
 Standard Markovian arguments (see Chapter 3 in \cite{PS}) indicate that $\mathbb{V}$ should be in $C^1$ and should satisfy the following requirements
\begin{align}
 \label{no2} \mathbb{V}_{xx}(x,y)=0 &&&  \mbox{if }(x,y)\in C,\,x<y,\\
 \label{no3} \mathbb{V}_y(x,x+)=0 &&& \mbox{for all }x,\\
 \label{no4} \mathbb{V}_x(x,y)=G_x(x,y) &&& \mbox{for all }(x,y)\in \partial C.
\end{align}
Note that equations \eqref{no2} and \eqref{no3} arise from the application of the generator of the Markov process $(X,Y)$ to the function $\mathbb{V}$; \eqref{no4} is the consequence of the principle of smooth fit. 
\smallskip

\emph{Step 3.} The key geometric properties of the continuation and stopping sets arises from the following arguments. First, by \eqref{0,0}, observe that
\begin{equation}\label{translation}
\begin{split}
\mathbb{V}(x,y)&=\sup_{\tau}\E_{0,0}\left[\left(x+\sup_{0\leq s\leq \tau}X_s\right)\vee y-(x+X_\tau)^2\right]\\
&=x-x^2+\sup_{\tau}\E_{0,0}\left[\left(\sup_{0\leq s\leq \tau}X_s\right)\vee (y-x)-X_\tau^2\right]\\
&=x-x^2+\mathbb{V}(0,y-x),
\end{split}
\end{equation}
where in the second line we have used the identity $\E_{0,0}X_\tau=0$. This yields the following translation property of $C$:  if $(x,y)\in C$ and $\lambda\geq -y$, then $(x+\lambda,y+\lambda)\in C$. Indeed, if $\mathbb{V}(x,y)>G(x,y)$, then
\begin{align*}
 \mathbb{V}(x+\lambda,y+\lambda)&=x+\lambda-(x+\lambda)^2+\mathbb{V}(0,y-x)\\
&=\lambda-2x\lambda-\lambda^2+\mathbb{V}(x,y)\\
&>\lambda-2x\lambda-\lambda^2+G(x,y)=G(x+\lambda,y+\lambda).
\end{align*}
By passing to the complement, $D$ enjoys the same translation property. 
The second observation is that if $(x,y)\in D$ and $y'>y$, then $(x,y')$ also lies in the stopping region. Indeed, if $a,\,b,\,c\in \R$ satisfy $b<c$, then $$a\vee b-b\geq a\vee c-c,$$ so for any stopping time $\tau$ we have the inequality
\begin{align*}
&\E_{0,0}\left[\left(x+\sup_{0\leq s\leq \tau}X_s\right)\vee y'-(x+X_\tau)^2\right]-G(x,y')\\
&\leq \E_{0,0}\left[\left(x+\sup_{0\leq s\leq \tau}X_s\right)\vee y-(x+X_\tau)^2\right]-G(x,y)\leq \mathbb{V}(x,y)-G(x,y)=0.
\end{align*}
Taking the supremum over all $\tau$, we obtain that $\mathbb{V}(x,y')\leq G(x,y')$, which implies that $(x,y')\in D$. Combining the above two observations, we see that there is a constant $a>0$ such that $$C=\{(x,y):0\leq y-x<a\}$$ and $$D=\{(x,y):y-x\geq a\}.$$ Note that the use of strict/non-strict inequalities comes from the fact that $C$ is open and $D$ is closed. This is due to the continuity of $\mathbb V$ and $G$.

\smallskip

\emph{Step 4.} Now, based on \eqref{no2}-\eqref{no4}, we provide the formula for the candidate for the value function, which will be denoted by $V$. By \eqref{no2} and \eqref{no4}, we see that if $(x,y)$ lies in the continuation set, then
$$ V(x,y)=G(y-a,y)+G_x(y-a,y)(x-y+a)=y+(y-a)^2-2(y-a)x.$$
Applying the condition in \eqref{no3} yields $a=\frac{1}{2}$. We have thus obtained that
$$ V(x,y)=\begin{cases}
y-x^2 & \mbox{if }y-x\geq \frac{1}{2},\\
y^2+\frac{1}{4}-(2y-1)x & \mbox{if }y-x<\frac{1}{2}.
\end{cases}$$

\smallskip

\emph{Step 5.} It is straightforward to see that the function $V$ obtained above is excessive (that is, it satisfies \eqref{no2}, \eqref{no3} and \eqref{no4}). Hence, by applying It\^o's formula, we have that $V\geq \mathbb V$. The reverse estimate is obtained by considering the stopping time given in \eqref{opt_stop}. This stopping time is integrable, even exponentially (see, for example, \cite{W}). Furthermore, for any $(x,y)\in \mathcal D$, the stopped process $(X^\tau,Y^\tau)$ evolves along the continuation region and It\^o's formula gives $$V(x,y)= \E_{x,y} V(X_{\tau},Y_{\tau})\leq \mathbb V(x,y).$$
This proves that $V=\mathbb V$. We see that the optimal stopping time \eqref{opt_stop} has the following interpretation: if the distance between $X$ and $Y$ is less than $\frac{1}{2}$, it is beneficial to wait; otherwise we should stop. This strategy makes intuitive sense as well: if $X$ is near its running maximum, then there is a high probability that the maximum will increase at any given moment (thus increasing the value $\mathbb{V}$), and the cost of waiting, expressed in terms of time or the increase of $\E X^2$, is relatively small. However, when the distance is large, it may take longer for $X$ to return to $Y$, so  the expected cost of waiting is too high and hence it is optimal to stop immediately.  These heuristics shall prove helpful in the sequel.

\subsection{On \eqref{optimal_stopping} for $n=1$} We now proceed with the study of the spider process. In the case $n=1$, the process coincides with the reflecting Brownian motion $|X|=(|X_t|)_{t\geq 0}$. Let $Y$ denote the corresponding two-sided maximal function $$Y=(Y_t)_{t\geq 0}=(\sup_{0\leq s\leq t}|X_s|)_{t\geq 0}.$$ Recalling \eqref{optimal_stopping} and invoking Wald's identity leads us to the optimal stopping problem
$$ \mathbb{U}=\sup \E (Y_\tau-X_{\tau}^2),$$
where the supremum is taken over all integrable stopping times $\tau$ of $|X|$. The analysis essentially proceeds along the same lines as in the previous subsection. Since the maximal function $Y$ above is two-sided, we modify the domain to $\mathcal{D}=\{(x,y): |x|\leq y\}$, introduce  the value function
$$ \mathbb{U}(x,y)=\sup \E_{x,y} (Y_\tau-|X_{\tau}|^2),\qquad (x,y)\in \mathcal{D},$$
and define the continuation and stopping regions $C$ and $D$ by the same formulas as before. Note that $G(x,y)=G(-x,y)$ and that the distribution of $(X,Y)$ under $\mathbb{P}_{x,y}$ is the same as that of $(-X,Y)$ under $\mathbb{P}_{-x,y}$. We thus conclude that
\begin{equation}\label{symmetry}
\mathbb{U}(x,y)=\mathbb{U}(-x,y)\qquad \mbox{for all }(x,y)\in \mathcal{D}.
\end{equation}
The key here difference is that the presence of two-sided maximal function disables equation \eqref{translation}, which had proved to be fundamental in the previous analysis.\\
\indent To get around this issue, we present the following reduction argument which shows that $\mathbb{U}$ and $\mathbb{V}$ coincide on a large part of the domain. Firstly, since $Y$ is not smaller than the one-sided maximal function, the direct comparison of the formulas for $\mathbb U$ and $\mathbb V$ gives that $\mathbb U\geq \mathbb V$ on $\mathcal{D}$. Next, suppose that $y_0$ is a nonnegative number such that $(0,y_0)\in D$. Repeating the reasoning from Step 3 above, we see that the entire half-line $$\ell=\{0\}\times [y_0,\infty)$$ is contained within $D$. Thus, for any $x\geq 0$, in the definition of $\mathbb{U}(x,y_0)$ one can restrict oneself to those stopping times $\tau$ for which the process $X^\tau$ does not go below zero: indeed, for other stopping times the process $(X^\tau,Y^\tau)$  crosses the line $\ell$, which is not optimal (see \eqref{opt_stop}). However, for such $\tau$, the process $Y$ coincides with the one-sided maximal function and hence by the very definition of $\mathbb U$ and $\mathbb V$ we have the desired reverse bound $\mathbb{U}(x,y_0)\leq \mathbb{V}(x,y_0)$. This in particular implies that $y_0\geq \frac{1}{2}$, since otherwise we would obtain 
$$ G(0,y_0)=\mathbb{U}(0,y_0)=\mathbb{V}(0,y_0)>G(0,y_0),$$
a contradiction. Denoting by $b$ the infimum of all $y_0$'s as above, we have that for $y \geq b$, $$\mathbb{U}(x,y)=\mathbb{V}(|x|,y).$$ \\
\indent We now apply Markovian arguments to obtain $\mathbb U_{xx}=0$ on $C$; one also may note that by the symmetry condition in \eqref{symmetry}, we have $\mathbb{U}_x(0,y)=0$ for $y<b$. These two observations imply that the function $\mathbb{U}$ must be constant on $\mathcal{D}\cap \{y<b\}$ and, by continuity, $\mathbb{V}$ must be constant on the line segment $[0,b]\times\{b\}$. This implies that $b=\frac{1}{2}$, leading us to the construction of the candidate function
$$ U(x,y)=\begin{cases}
y-x^2 & \mbox{if }y\geq x+ \frac{1}{2},\\
y^2+\frac{1}{4}-(2y-1)x & \mbox{if }x+\frac{1}{2}> y\geq \frac{1}{2},\\
\frac{1}{2} & \mbox{if }y< \frac{1}{2}.
\end{cases}$$
It is straightforward to check that $U$ is excessive and hence that $U\geq \mathbb{U}$; the reverse bound is obtained by considering the stopping time in \eqref{opt_stop}. The optimal strategy is to wait until the distance between $|X|$ and its running maximum is at least $\frac12$. Intuitively, this remains perfectly consistent with the strategy for the previous problem.

\subsection{The case $n=2$} Here the analysis will be a more involved, but again the special function $\mathbb{V}$ will play a prominent role. 

\smallskip

\emph{Step 1.} The spider process $S$ coincides with the standard one-dimensional Brownian motion, which we shall denote again by $X$. Let  $Y=(Y_t)_{t\geq 0}$ and $Z=(Z_{t})_{t\geq 0}$ be the running maximum and the running infimum of $X$; that is, for $t \geq 0$, $$Y_t=\sup_{0\leq s\leq t}X_s$$ and $$Z_{t}=\inf_{0\leq s\leq t}X_s.$$ Motivated by \eqref{optimal_stopping} and Wald's identity, we introduce the optimal stopping problem  
$$ \mathbb{U}=\sup \E (Y_\tau-Z_{\tau}-X_{\tau}^2),$$
where the supremum is taken over all integrable stopping times $\tau$ of $X$. It is important to note that, in contrast to the previous considerations, there are now three variables involved. 

To apply the general theory of optimal stopping, we extend the triple $(X,Y,Z)$ to a Markov family on the state space $$\mathcal{D}=\{(x,y,z): 0,\,x\in [z,y]\}.$$ Let the corresponding family of initial distributions be denoted by $(\mathbb{P}_{x,y,z})_{(x,y,z)\in \mathcal{D}}.$ Having done that, we introduce the gain function $G(x,y,z)=y-z-x^2$ and the value function
\begin{equation}\label{stop_n=2'}
 \mathbb{U}(x,y,z)=\sup \E_{x,y,z}G(X_\tau,Y_\tau,Z_{\tau}).
\end{equation}
Here the supremum is taken over all $\mathbb{P}_{x,y,z}$-integrable stopping times $\tau$ of $X$. 
The associated continuation and the instantaneous stopping regions are given by 
\begin{equation}\label{CD}
\begin{split}
 C&=\big\{(x,y,z)\in \mathcal{D}\,:\, \mathbb{U}(x,y,z)>G(x,y,z)\big\},\\
 D&=\big\{(x,y,z)\in \mathcal{D}\,:\, \mathbb{U}(x,y,z)=G(x,y,z)\big\},
\end{split}
\end{equation}
and the optimal stopping time in \eqref{stop_n=2'} is
\begin{equation}\label{opt_stop2}
\tau=\inf\{t\geq 0:(X_t,Y_t,Z_t)\in D\}.
\end{equation}

\smallskip

\emph{Step 2.} We now provide an initial comparison of the functions $\mathbb U$ and $\mathbb V$, exploiting a similar argument as in the case $n=1$. We begin with the observation that both the inequalities $Y_\tau\geq y$ and $Z_\tau\leq z$ hold $\mathbb{P}_{x,y,z}$-almost surely. This implies that
\begin{equation} \label{estimate1}
 \E_{x,y,z}(Y_\tau-Z_\tau-X_\tau^2)\geq \E_{x,y,z}(Y_\tau-X_\tau^2)-z
 \end{equation}
 and
 \begin{equation} \label{estimate2}
 \E_{x,y,z}(Y_\tau-Z_\tau-X_\tau^2)\geq \E_{x,y,z}(-Z_\tau-X_\tau^2)+y.
 \end{equation}
 By the definition of $\mathbb U$ and $\mathbb V$,  \eqref{estimate1} gives
\begin{equation}\label{U>V}
\mathbb U(x,y,z)\geq \mathbb V(x,y)-z.
\end{equation}
To see the consequence of \eqref{estimate2}, note that, under $\mathbb{P}_{x',y,z}$, $(X,-Z)$ has the same distribution as $(X,Y)$ under $\mathbb{P}_{-x,-z,-y}$. We therefore obtain
\begin{equation}\label{U>>V}
\mathbb U(x,y,z)\geq \mathbb V(x,-z)+y.
\end{equation}
\indent Indeed, both \eqref{U>V} and \eqref{U>>V} can be reversed on a large part of the domain. Let $z<0<y$ be fixed numbers and suppose that there is $x\in (z,y)$ such that the state $(x,y,z)$ belongs to the stopping domain. Then the whole half-line $$\ell=\{x\}\times [y,\infty)\times \{z\}$$ is entirely contained within $D$ (repeating the argument from Step 3 in Subsection 2.1). Next, suppose that $x'>x$. Then in the definition of $\mathbb{U}(x',y,z)$ it is enough to consider only those stopping times $\tau$ for which the process $(X^\tau)_{t\geq 0}$ does not go below $x$. Indeed, for other stopping times the process $(X^\tau,Y^\tau,Z^\tau)$ crosses the line $\ell$, which is not optimal (see \eqref{opt_stop2}). However, for such $\tau$ the running infimum $Z^\tau$ will not change, so
$$ \mathbb{U}(x',y,z)=\sup \E_{x',y,z}(Y_\tau-X_\tau^2)-z\leq \mathbb V(x',y)-z.$$
An analogous argument works for $x'<x$: in this case, when studying \eqref{stop_n=2'}, one may restrict oneself to stopping times $\tau$ for which $X^\tau$ does not go above $x$, which keeps $Y^\tau$ fixed and yields the desired reverse identity
$$ \mathbb{U}(x',y,z)=\sup \E_{x',y,z}(-Z_\tau-X_\tau^2)+y\leq \mathbb V(x',-z)+y.$$

\smallskip

\emph{Step 3.} Note that \eqref{U>V} and \eqref{U>>V} imply that if $y-z<1$, then $(x,y,z)\in C$ for all $x\in [z,y]$. Indeed, for any such $x$ we have $x-z<\frac{1}{2}$ or $y-x<\frac{1}{2}$, and hence $\mathbb V(x,y)>y-x^2$ or $\mathbb V(x,-z)\geq -z-x^2$; this gives 
$$\mathbb U(x,y,z)> G(x,y,z).$$ It is useful to note that this in perfect consistence with the optimal strategies described at the end of the previous two subsections: if $y-z<1$, then the distance between $x$ and $y$ or the distance between $x$ and $z$ is less than $\frac12$, and hence it is beneficial to wait. Actually, this observation also suggests what to do if $y-z\geq 1$. If both $x-z$ and $y-x$ are at least $\frac12$, one should stop; otherwise, wait. In other words, by the analysis carried out in the previous step, we obtain that the candidate $U$ for the value function satisfies, if $y-z\geq 1$,
$$ \mathbb{U}(x,y,z)=\begin{cases}
\mathbb V(x,y)-z & \mbox{if }y-x<x-z,\\
\mathbb V(x,-z)+y & \mbox{if }y-x\geq x-z.
\end{cases}$$
For $y-z<1$, one exploits Markovian arguments and obtains the system of equations
\begin{align*}
  U_{xx}(x,y,z)=0 &&&  \mbox{if }z<x<y,\\
  U_y(x,x+,z)=0 &&& \mbox{for all }z<0<x,\\
  U_z(x,y,x-)=0 &&& \mbox{for all }x<0<y.
\end{align*}
This system can be solved explicitly (cf. \cite{DGM,E}, see also Section 3 below): we obtain
$$ U(x,y,z)=y-z-x(y+z)+\frac{(y-1)^2+(z+1)^2}{2}-\frac{1}{4}.$$

\smallskip

\emph{Step 4.} The analysis is completed by showing that $U=\mathbb U$. This is done as we have done so previously: one checks that $U$ is excessive and hence $U\geq \mathbb U$. The reverse bound follows from the construction, since $U$ is obtained by exercising the optimal strategy described above. We omit the details, referring the interested reader to \cite{DGM,E}.

\section{On the search for the value function for $n\geq 3$}

Equipped with the above machinery and intuition, we proceed to the analysis of the case $n=3$. The purpose of this section is to obtain a candidate $U$ for the value function associated with the appropriate optimal stopping problem. The reasoning rests on a number of guesses and assumptions that may (at least at first glance) seem imprecise. However, the reader should keep in mind that our explicit purpose here is \textit{only to guess} an appropriate special function. The necessary rigorous analysis will be presented in Section 4. Again, for purposes of clarity, we split the reasoning into intermediate steps.

\smallskip

\emph{Step 1.} First, we need to specify the underlying Markov process which will be subject to the optimal stopping procedure. Of course, we could consider the process $$(X,Y^{(1)},Y^{(2)},\ldots,Y^{(n)}),$$ where $X$ takes the values in $$R_1\cup R_2\cup \ldots \cup R_n,$$ and  $Y^{(j)}_t$ measures the length of $j$-th rib up to time $t$, but this process has a very complicated structure. Fortunately, there is an alternative choice for which the state space is simpler: a three-dimensional structure. As a starting point, note that the diameter of the spider process depends only on the behavior of two longest ribs and hence, as it was for $n=2$, it is natural to try to find some representation of $S$ on the real line. For $t\geq 0$, distinguish the longest rib
$$ Y_t=\max_{1\leq j\leq n}|T_t(\omega)\cap R_j|$$
(where $T_t(\omega)$ was defined in \eqref{trajectory}) and let $-Z_t$ (note the minus sign) be the corresponding second longest rib, so that $D_t=Y_t-Z_t$. Now, to  define $X$, we first set $|X_t|$ to be the distance of the spider process $S_t$ from the origin. Furthermore, if $S_t$ belongs to the ``running longest rib'', we assume that $X_t\geq 0$; otherwise, we assume that $X$ is negative. In other words, we copy the running longest rib on the positive half-line, while all the remaining ribs are glued together and copied on $(-\infty,0]$. For a graphical depiction of the above arguments, see Figure \ref{fig22} below.

The process $X$ can be interpreted in the language of skew Brownian motion (see, for example, \cite{Shepp}). Given $\alpha\in [0,1]$, the  $\alpha$-skew Wiener process can be obtained from reflecting Brownian motion by changing (independently) the sign of each excursion with probability $\alpha$. Thus, $0$-skew Wiener process is reflecting Brownian motion, while $\frac{1}{2}$-skew Wiener process is the usual Brownian motion. The $\alpha$-skew Wiener process behaves like a usual Wiener process except for the asymmetry at the origin: if located at zero, then for any $s>0$ the process has probability $\alpha$ of reaching $-s$ before $s$. 

Note that the process $X$ defined above is a $1-n^{-1}$-skew Brownian motion, which possesses the additional jump part: if for a given $t>0$, its left limit $X_{t-}$ equals $-Y_t$, then $X_t$ changes its sign, moving to $Y_t$. This discontinuity (or ``phase-transition'') corresponds to the scenario in which the second longest rib becomes the longest.
\begin{figure}[htbp]
\begin{center}
\includegraphics[scale=0.29]{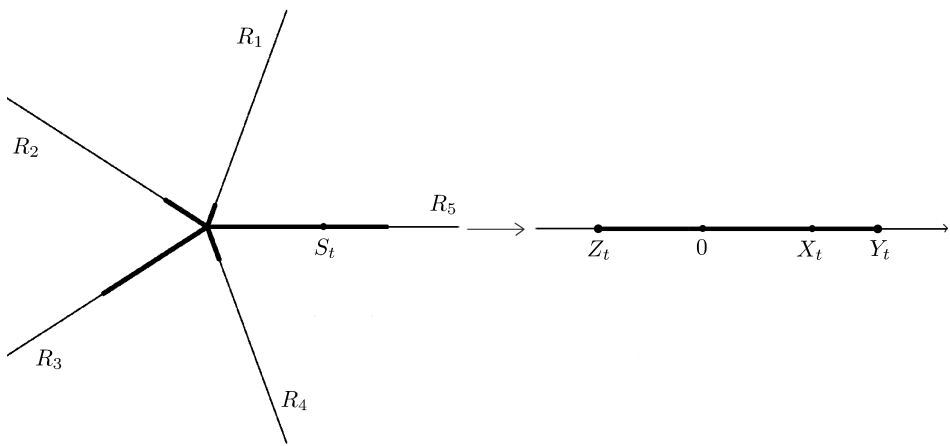}
\caption{The spider process ($n=5$) and its transformation to the skew Brownian motion $X$ with jumps. The ray $R_5$ containing the longest rib has been copied onto the positive half-line; the remaining rays $R_1$-$R_4$ have been glued together and copied onto the negative half-line. If $X$ reaches $-Y_t$ before $Y_t$, it then jumps to $Y_t$. }\label{fig22}
\end{center}
\end{figure}

\smallskip

\emph{Step 2.} We now gather some basic information about the behavior of the triple $(X,Y,Z)$. It is straightforward to check that this is a time-homogeneous, right-continuous strong Markov process on the state space $$\mathcal{D}=\{(x,y,z): z\leq 0\leq y,\,z\leq x\leq y,\,y+z\geq 0\}.$$ As usual, we will denote by $(\mathbb{P}_{x,y,z})_{(x,y,z)\in \mathcal{D}}$ the corresponding family of initial distributions such that 
$$\mathbb{P}_{x,y,z}((X_0,Y_0,Z_0)=(x,y,z))=1.$$
It is immediate that the process $(X,Y,Z)$ enjoys the following Brownian scaling property.

\begin{lemma}
For any $\lambda>0$, the process $$t\mapsto (\lambda X_{t\lambda^{-1/2}},\lambda  Y_{t\lambda^{-1/2}},\lambda Z_{t\lambda^{-1/2}})$$ has the same law under $\mathbb{P}_{x,y,z}$ as does $(X,Y,Z)$ under $\mathbb{P}_{\lambda x,\lambda y,\lambda z}$.
\end{lemma}

We will also need the following property.

\begin{lemma}\label{Property XYZ}
Let $y<1/2$ and $\sigma=\inf\{t>0:Y_t\geq 1/2\}$. Then the distribution of $Z_\sigma$ under $\mathbb{P}_{x,y,z}$ is determined by
$$ \mathbb{P}_{x,y,z}(Z_\sigma\leq s)=\begin{cases}
0 & \mbox{if }s<-\frac{1}{2},\\
\displaystyle \frac{(n-1)(1-2x)}{n-1-2s} & \mbox{if }x\geq 0,\,-y\leq s<z,\medskip \\
\displaystyle \frac{n-1-2x}{n-1-2s} & \mbox{if }x< 0,\,-y\leq s<z,\medskip\\
\displaystyle \frac{(n-1)(1+2s)}{n-1-2s} & \mbox{if }-\frac{1}{2}\leq s<-y,\\
1 & \mbox{if }s\geq z.
\end{cases}$$
\end{lemma}
\begin{proof}
It suffices to prove the formula for $-\frac{1}{2}\leq s<z$, since for the remaining $s$ the claim is obvious. We consider three separate cases. 

\smallskip

\noindent \underline{Case I}: Suppose that $x\geq 0$ and $y\geq -s$. By the law of total probability, we obtain
\begin{equation}\label{total}
 \mathbb{P}_{x,y,z}(Z_\sigma> s)=2x+(1-2x)\mathbb{P}_{0,y,z}(Z_\sigma > s).
\end{equation}
The above equation contains two disjoint scenarios. The process $X$ may get to $1/2$ before it visits $0$: this occurs with probability $2x$ and automatically implies that $Z_\sigma>s$. The second possibility is that $X$ drops to zero before it reaches $1/2$. Then,  no matter how much $Y$ increased, the set $\{Z_\sigma>s\}$ has conditional probability $\mathbb{P}_{0,y,z}(Z_\sigma > s)$; indeed, the latter does not depend on the value of $y\in [-s,1/2)$. 
To compute this probability, note that since $y\geq -s$, we have 
\begin{equation}\label{total2}
\mathbb{P}_{0,y,z}(Z_\sigma > s)=\mathbb{P}_{-s,y,z}(Z_\sigma > s)/n.
\end{equation}
Indeed, the inequality $Z_\sigma>s$ means that when $X$ reaches $-s$, the spider process is on the longest rib: by symmetry, the probability of this scenario is $1/n$. After that, no matter how much $Z$ has dropped, the event $\{Z_\sigma>s\}$ occurs with the conditional probability equal to $\mathbb{P}_{-s,y,z}(Z_\sigma > s)$. Now, applying \eqref{total} with $x=-s$ and combining it with \eqref{total}, we obtain that $$ \mathbb{P}_{0,y,z}(Z_\sigma>s)=-2s/(n-1-2s),$$ or
$$ \mathbb{P}_{0,y,z}(Z_\sigma\leq s)=\frac{n-1}{n-1-2s}.$$
Plugging this into \eqref{total} yields
$$ \mathbb{P}_{x,y,z}(Z_\sigma\leq s)=\frac{(n-1)(1-2x)}{n-1-2s}.$$

\smallskip

\noindent \underline{Case II}:  Next, assume that $x<0$ and $y\geq -s$. The inequality $Z_\sigma>s$ implies that $X$ must rise to $0$ before it drops to $s$, and the change in $Z$ is irrelevant, so
$$ \mathbb{P}_{x,y,z}(Z_\sigma>s)=\left(1-\frac{x}{s}\right)\mathbb{P}_{0,y,z}(Z_\sigma>s)=\frac{2(x-s)}{n-1-2s}.$$

\smallskip

\noindent \underline{Case III}:  Finally, suppose that $y<-s$. Then conditioning on the moment at which $X$ visits $-s$ for the first time, we obtain
$$\mathbb{P}_{x,y,z}(Z_\sigma>s)=\mathbb{P}_{-s,-s,z}(Z_\sigma>s).$$ Again, the drop in $Z$ is not important and we may write $z$ in the lower index on the right. Hence
$$ \mathbb{P}_{x,y,z}(Z_\sigma\leq s)=\mathbb{P}_{-s,-s,z}(Z_\sigma\leq s)=\frac{(n-1)(1+2s)}{n-1-2s},$$
where the latter equality follows from the analysis in Case I.
\end{proof}

\emph{Step 3.} We proceed to the study of \eqref{optimal_stopping}. As before, we extend it to an arbitrary starting point $(x,y,z)\in \mathcal{D}$, setting
\begin{equation}\label{stop_n}
 \mathbb{U}(x,y,z)=\sup \E_{x,y,z}G(X_\tau,Y_\tau,Z_\tau),
\end{equation}
where $$G(x,y,z)=y-z-x^2$$ and the supremum is taken over all integrable stopping times $\tau$ of $X$. 
Markovian arguments show that $\mathbb U$ satisfies the following system of equations
\begin{align}
\label{21} \mathbb{U}_{xx}(x,y,z)&=0 && \mbox{if }(x,y,z)\in C,\quad z<x<y,\,x\neq 0,\\
\label{22} \mathbb{U}_y(y,y+,z)&=0 && \mbox{for all }y>0,\\
\label{23} \mathbb{U}_z(z,y,z-)&=0 && \mbox{for all }z<0,\\
\label{24} (n-1)\mathbb U_x(0-,y,z)&=\mathbb U_x(0+,y,z) && \mbox{if }(0,y,z)\in C.
\end{align}
As the direct consequence of \eqref{21}, the stopping set has the property that if it contains two points of the form $(x,y,z)$ and $(x',y,z)$, then it automatically contains the entire line segment which joins these two points. Otherwise, by the concavity of the function $x\mapsto G(x,y,z)$, this would violate the inequality $\mathbb{U}\geq G$.

\smallskip

\emph{Step 4.} Our construction for the candidate $U$ for the value function will be based on the guess of the optimal stopping strategy. Equipped with the analysis in the case $n=2$, a naive idea is to try to proceed analogously, i.e., consider the optimal stopping times $\tau$ which consists of two stages: 

\smallskip

Stage 1. Wait until the difference between $Y$ and $Z$ is equal to $1$;

Stage 2. Wait until $Y-X$ and $X-Z$ are both larger than $\frac{1}{2}$.

\smallskip

Some thought reveals that this cannot be the optimal strategy. To see this, suppose that the first stage is over and then, after some time, we have $X=0$ and $Z\in (-\frac{1}{2},0)$. Because of the asymmetry of the skew Brownian at zero (which in our case ``pushes'' the process on the negative side),  the cost of waiting for $X$ to reach $Z$ is \emph{lower} than in the symmetric case, so the margin $\frac{1}{2}$ should be increased, at least if at the end of Stage 1 we have $X=Z$. 

On the other hand, it is natural to expect that the above strategy is not far from optimal. It seems plausible to try the following general two-step procedure:

\smallskip

Stage 1. Wait until $Y$ and $Z$ become ``distant'';

Stage 2. Wait until $f(Y,Z)\leq X\leq g(Y,Z)$ for some functions $f$ and $g$.

\smallskip

\noindent Note that in the light of the above arguments, we must have $Y-Z>1$ and hence $Y>\frac{1}{2}$ at the end of the first stage (we have $Y\geq -Z$ almost surely).

\smallskip

\emph{Step 5.} We now turn to the study of some basic properties of $f$ and $g$. First, note that $f$ depends only on $z$ and $g$ depends only on $y$.  The idea behind this is as follows. Suppose that  $(x,y,z)\in D$; then $(x,y',z)$ and $(x,y,z')$ also lie in the stopping set, provided $y'>y$ and $z'<z$ (the argument is the same as in the case $n=2$). So, when computing $\mathbb{U}({z,y,z})$ we may restrict ourselves to those stopping times $\tau$ for which $X$ does not cross $f(y,z)$; for such $\tau$ the process $Y^\tau$ is constant and hence for $x<f(y,z)$ we have $$\mathbb{U}(x,y,z)=\sup \E_{x,y,z}(-Z_\tau-X_\tau^2)+y.$$ The problem thus reduces to the optimal stopping of $X$ and $Z$: thus the stopping boundary cannot depend on $y$ and hence $f(Y,Z)=f(Z)$. We can now go one step further: if $f(z)\leq 0$, then for $\tau$ as above $-Z^\tau$ is the one-sided maximal function of $-X^\tau$ and hence $$\mathbb{U}(x,y,z)=\mathbb V(-x,-z)+y,$$ so in particular $f(z)=z+\frac{1}{2}$.  
A similar argument shows that $g(Y,Z)=g(Y)$; however, since $y> \frac{1}{2}$ (see the end of the previous step), we obtain $g(y)=y-\frac{1}{2}$ for all $y$.

\smallskip


\emph{Step 6.} Now we will find the formula for $f$ for $z$ close to zero (so that $f(z)>0$). To this end, we will show that $f$ satisfies an appropriate ordinary differential equation. We thus fix such a $z$. By \eqref{21}, \eqref{24} and the principle of smooth fit 
$$\mathbb U_x(f(z)-,y,z)=G_x(f(z)+,y,z)=-2f(z),$$ 
we obtain the identity 
$$ U(x,y,z)=\begin{cases}
y-z+f(z)^2-2f(z)x &\mbox{ for }x\in (0,f(z)),\\
\displaystyle y-z+f(z)^2-\frac{2f(z)x}{n-1} &\mbox{ for }x\in (z,0).
\end{cases}$$
Applying \eqref{23}, we get that
$$ 2f'(z)\left[f(z)-\frac{z}{n-1}\right]=1.$$
Note the initial condition $f(-\frac{1}{2})=0$, which comes from the case $z\leq -\frac{1}{2}$ considered above. 
This differential equation can be easily solved: the substitution $z=\varphi(s)=f^{-1}(s)$ transforms it into the linear equation
$$ 2\left(s-\frac{\varphi(s)}{n-1}\right)=\varphi'(s),\qquad \varphi(0)=-\frac{1}{2},$$
whose explicit solution is
$$ \varphi(s)=(n-1)s-\frac{(n-1)^2}{2}+\frac{n(n-2)}{2}\exp\left(-\frac{2s}{n-1}\right).$$
Hence, for $z>-\frac{1}{2}$, $f$ is the inverse to the above function.

\smallskip

\emph{Step 7.} We are now ready to guess the final form of the optimal strategy. We have already constructed an appropriate lower and upper boundary functions $f$ and $g$. Taking the above discussion into account, we formulate the procedure as follows.

\smallskip

Stage 1. Wait until the equality $f(Z)\leq g(Y)$ is observed for the first time.

Stage 2. Wait until $f(Z)\leq X\leq g(Y)$.

\smallskip

The remaining part of the analysis is devoted to the explicit evaluation of the value function associated with this strategy. In other words, we shall henceforth set
$$ U(x,y,z)=\E_{x,y,z}G(X_\tau,Y_\tau,Z_\tau),$$
where $\tau$ is given as the combination of Stage 1 and Stage 2 above. The discussion we have already carried out gives the following (partial) formula for $U$.

\begin{cor}\label{partU}
If $y\geq 1/2$ and $z\leq \varphi(y-\frac{1}{2})$, then 
\begin{equation}\label{form1}
 U(x,y,z)=\begin{cases}
y-z-2\left(z+\frac{1}{2}\right)x+\left(z+\frac{1}{2}\right)^2 & \mbox{if }x\leq \varphi^{-1}(z)\leq 0,\\
y-z-\frac{2\varphi^{-1}(z)}{n-1}x+\big(\varphi^{-1}(z)\big)^2 & \mbox{if }x\leq 0\leq \varphi^{-1}(z),\\
y-z-{2\varphi^{-1}(z)}x+\big(\varphi^{-1}(z)\big)^2 & \mbox{if }0\leq x\leq \varphi^{-1}(z),\\
y-z-x^2 & \mbox{if }\varphi^{-1}(z)\leq x\leq y-\frac{1}{2},\\
y-z-2\left(y-\frac{1}{2}\right)x+\left(y-\frac{1}{2}\right)^2 & \mbox{if }x\geq y-\frac{1}{2}.
\end{cases}
\end{equation}
\end{cor}

It remains to find the formula for $U$  for $z>\varphi(y-\frac{1}{2})$. We consider the cases $y\geq \frac{1}{2}$ and $y<\frac{1}{2}$ separately.

\smallskip

\emph{Step 8.} First we study the case $y\geq \frac{1}{2}$; this is the most elaborate part. We begin with a formula for $U_y$.

\begin{lemma}\label{lem34}
Let $y\geq 1/2$ and $z>\varphi(y-\frac{1}{2})$. The function $U$ satisfies 
\begin{equation}\label{uy}
 U_y(x,y,z)=\begin{cases}
\displaystyle \frac{(n-1)(y-x)}{(n-1)y-\varphi(y-\mbox{$\frac{1}{2}$})} & \mbox{if }x\geq 0, \smallskip\\
\displaystyle \frac{(n-1)y-x}{(n-1)y-\varphi(y-\mbox{$\frac{1}{2}$})} & \mbox{if }x\leq 0.
\end{cases}
\end{equation}
\end{lemma}
\begin{proof}
It suffices to prove the formula for $x\geq 0$: indeed, by Markovian arguments we see that $U$ satisfies \eqref{21} and \eqref{24} (with $\mathbb U$ replaced by $U$), so 
$$U(x,y,z)=\begin{cases}
U(0,y,z)+U_x(0-,y,z)x & \mbox{if }x\leq 0,\\
U(0,y,z)+(n-1)U_x(0-,y,z)x & \mbox{if }x\geq 0
\end{cases}$$
and
$$U_y(x,y,z)=\begin{cases}
U_y(0,y,z)+U_{xy}(0-,y,z)x & \mbox{if }x\leq 0,\\
U_y(0,y,z)+(n-1)U_{xy}(0-,y,z)x & \mbox{if }x\geq 0.
\end{cases}$$
Hence if \eqref{uy} is valid for $x\geq 0$, it automatically holds for $x<0$ as well. 

Therefore, we shall henceforth assume that $x\geq 0$. Our plan is to write 
\begin{equation}\label{differential}
 U_y(x,y+,z)=\lim_{\delta \downarrow 0}\frac{U(x,y+\delta,z)-U(x,y,z)}{\delta}
\end{equation}
and to analyze the expectations defining $U(x,y,z)$ and $U(x,y+\delta,z)$. To this end, we fix a small $\delta>0$ (so that $z\geq \varphi(y+\delta-\frac{1}{2})$) and consider the events
\begin{align*}
A_1&=\{\mbox{the trajectory of $X$ reaches $y+\delta$ before $\varphi(y+\delta-1/2)$}\},\\
A_2&=\{\mbox{the trajectory of $X$ reaches $\varphi(y+\delta-1/2)$ before $y$}\},\\
A_3&=\big\{\mbox{the trajectory of $X$ reaches $y$ before $\varphi(y+\delta-1/2)$},\\
&\qquad   \mbox{but after that it reaches $\varphi(y+\delta-1/2)$ before $y+\delta$}\big\}.
\end{align*}
Of course $A_1$, $A_2$, $A_3$ are pairwise disjoint and their union has probability $1$. Now, we write
\begin{align*}
U(x,y+\delta,z)&=\E_{x,y+\delta,z}(Y_\tau-Z_\tau-X_\tau^2)=I_1+I_2,
\end{align*}
where
\begin{align*}
I_1&=\E_{x,y+\delta,z}\Big((Y_\tau-Z_\tau-X_\tau^2)1_{A_1}\Big),\qquad I_2=\E_{x,y+\delta,z}\Big((Y_\tau-Z_\tau-X_\tau^2)1_{A_1^c}\Big),
\end{align*}
and $A_1^c=\Omega\setminus A_1$ is the complement of $A_1$. By the Markov property, we see that 
\begin{equation}\label{chain}
\begin{split}
I_2&=\E_{x,y+\delta,z}\Big((Y_\tau-Z_\tau-X_\tau^2)|A_1^c\Big)\mathbb{P}_{x,y+\delta,z}(A_1^c)\\
&=U\big(\varphi(y+\delta-\mbox{$\frac{1}{2}$}),y+\delta,\varphi(y+\delta-\mbox{$\frac{1}{2}$})\big)\cdot \frac{(n-1)(y+\delta-x)}{(n-1)(y+\delta)-\varphi(y+\delta-\frac{1}{2})}.
\end{split}
\end{equation}
We write down a similar splitting for $U$ as $U(x,y,z)=J_1+J_2+J_3$, with
\begin{align*}
J_k= \E_{x,y,z}\Big((Y_\tau-Z_\tau-X_\tau^2)1_{A_k}\Big),\qquad k=1,\,2,\,3.
\end{align*}
A crucial observation is that
\begin{align*}
 I_1=\E_{x,y+\delta,z}\Big((Y_\tau-Z_\tau-X_\tau^2)1_{A_1}\Big)&=\E_{x,y,z}\Big((Y_\tau \vee(y+\delta)-Z_\tau-X_\tau^2)1_{A_1}\Big)\\
&=\E_{x,y,z}\Big((Y_\tau-Z_\tau-X_\tau^2)1_{A_1}\Big)=J_1,
\end{align*}
since $Y_\tau\geq y+\delta$ on $A_1$ (by the very definition of this event). Furthermore, arguing as in \eqref{chain}, we obtain
\begin{align*}
J_2&=U\big(\varphi(y+\delta-\mbox{$\frac{1}{2}$}),y,\varphi(y+\delta-\mbox{$\frac{1}{2}$})\big)\cdot \frac{(n-1)(y-x)}{(n-1)y-\varphi(y+\delta-\frac{1}{2})}.
\end{align*}
Finally, to work with $J_3$, we rewrite $A_3$ as the intersection of the following two events
\begin{align*}
A_3^1&=\big\{\mbox{the trajectory of }X\mbox{ reaches $y$ before $\varphi(y+\delta-1/2)$}\big\},\\
 A_3^2&=\big\{\mbox{having visited $y$, the trajectory of $X$ reaches $\varphi(y+\delta-\mbox{$\frac{1}{2}$})$ before $y+\delta$}\big\}.
\end{align*}
Then, using the Markov property, we compute that
\begin{align*}
&\E_{x,y,z}\Big((Y_\tau-Z_\tau-X_\tau^2)1_{A_3}\Big)=\E_{x,y,\varphi(y+\delta-1/2)}\Big((Y_\tau-Z_\tau-X_\tau^2)1_{A_3^1\cap A_3^2}\Big)\\
&=\E_{y,y,\varphi(y+\delta-1/2)}\Big((Y_\tau-Z_\tau-X_\tau^2)1_{A_3^2}\Big)\mathbb{P}_{x,y,z}(A_3^1)\\
&=\E_{y,y,\varphi(y+\delta-1/2)}\Big((Y_\tau-Z_\tau-X_\tau^2)1_{A_3^2}\Big)\left[1-\frac{(n-1)(y-x)}{(n-1)y-\varphi(y+\delta-1/2)}\right].
\end{align*}
To analyze the latter expectation, note that on the set $A_3^2$, when $X$ gets to $\varphi(y+\delta-\frac{1}{2})$, then the value of $Y$ lies between $y$ and $y+\delta$. Consequently, 
\begin{align*}
&\E_{y,y,\varphi(y+\delta-1/2)}\Big((Y_\tau-Z_\tau-X_\tau^2)1_{A_3^2}\Big)\\
&=U\big(\varphi(y+\delta-\mbox{$\frac{1}{2}$}),y,\varphi(y+\delta-\mbox{$\frac{1}{2}$})\big)\cdot \frac{(n-1)\delta}{(n-1)(y+\delta)-\varphi(y+\delta-\frac{1}{2})}+o(\delta).
\end{align*}
Plugging all the above into \eqref{differential}, we obtain
$$ U_y(x,y+,z)=\lim_{\delta \downarrow 0}\frac{I_1+I_2-(J_1+J_2+J_3)}{\delta}=\lim_{\delta\downarrow 0} \left(\frac{I_2-J_2}{\delta}-\frac{J_3}{\delta}\right).$$
We also have the identity
\begin{align*}
 \lim_{\delta\downarrow 0}\frac{I_2-J_2}{\delta}&=\frac{\partial}{\partial w}\left[U\big(\varphi(y-\mbox{$\frac{1}{2}$}),w,\varphi(y-\mbox{$\frac{1}{2}$})\big)\frac{(n-1)(w-x)}{(n-1)w-\varphi(y-\mbox{$\frac{1}{2}$})}\right]\Bigg|_{w=y}\\
 &=U_y\big(\varphi(y-\mbox{$\frac{1}{2}$}),y,\varphi(y-\mbox{$\frac{1}{2}$})\big)\frac{(n-1)(y-x)}{(n-1)y-\varphi(y-\mbox{$\frac{1}{2}$})}\\
 &\quad +U\big(\varphi(y-\mbox{$\frac{1}{2}$}),y,\varphi(y-\mbox{$\frac{1}{2}$})\big)\frac{\partial}{\partial w}\left[\frac{(n-1)(w-x)}{(n-1)w-\varphi(y-\mbox{$\frac{1}{2}$})}\right]\Bigg|_{w=y}.
\end{align*}
We have already computed above (see \eqref{form1}) that $$U_y\big(\varphi(y-\mbox{$\frac{1}{2}$}),y,\varphi(y-\mbox{$\frac{1}{2}$})\big)=1.$$ Furthermore, it is straightforward to check that the term
$$ U\big(\varphi(y-\mbox{$\frac{1}{2}$}),y,\varphi(y-\mbox{$\frac{1}{2}$})\big)\frac{\partial}{\partial w}\left[\frac{(n-1)(w-x)}{(n-1)w-\varphi(y-\mbox{$\frac{1}{2}$})}\right]\Bigg|_{w=y}$$ 
is precisely $ \lim_{\delta\downarrow 0}{J_3}/{\delta}$. Combining all of the above arguments,  we obtain the desired claim. This concludes the proof.
\end{proof}

Lemma \ref{lem34} allows us to extend the formula for $U$ to the domain $$\{(x,y,z):y\geq 1/2,\,z>\varphi(y-\frac{1}{2})\}.$$

\begin{cor}\label{cor+int}
If $y\geq 1/2$ and $z>\varphi(y-\frac{1}{2})$, then
$$ U(x,y,z)=\begin{cases}
\displaystyle U(x,f(z)+\mbox{$\frac{1}{2}$},z)-\int_y^{f(z)+1/2}\frac{(n-1)(s-x)\mbox{d}s}{(n-1)s-\varphi(y-\frac{1}{2})} &\mbox{if }x\geq 0,\smallskip\\
\displaystyle U(x,f(z)+\mbox{$\frac{1}{2}$},z)-\int_y^{f(z)+1/2}\frac{((n-1)s-x)\mbox{d}s}{(n-1)s-\varphi(y-\frac{1}{2})} &\mbox{if }x\leq 0.
\end{cases}$$
\end{cor}

\emph{Step 9.} This is the final part, concerning the case $y<1/2$, and it is much simpler. If we denote $\sigma=\inf\{t\geq 0:Y_t=\frac{1}{2}\}$, then the Markov property gives
$$ {U}(x,y,z)=\mathbb{E}_{x,y,z}U(X_\sigma,Y_\sigma,Z_\sigma)=\mathbb{E}_{x,y,z}U(\mbox{$\frac{1}{2},\frac{1}{2}$},Z_\sigma).$$
To compute the latter expectation, we apply Lemma \ref{Property XYZ} and immediately obtain the following.

\begin{cor}\label{const2}
If $x<0$ and $y<1/2$, then we have
\begin{equation}\label{Usmall1}
\begin{split}
 U(x,y,z)&=U(\mbox{$\frac{1}{2}$,$\frac{1}{2}$},-y)\cdot \frac{2((n-1)y-x)}{n-1+2y}+U(\mbox{$\frac{1}{2}$,$\frac{1}{2}$},z)\cdot \frac{2(x-z)}{n-1-2z}\\
&+\int_{-1/2}^{-y} U(\mbox{$\frac{1}{2}$,$\frac{1}{2}$},s)\cdot \frac{2n(n-1)}{(n-1-2s)^2}\mbox{d}s+\int_{-y}^{z} U(\mbox{$\frac{1}{2}$,$\frac{1}{2}$},s)\cdot \frac{2(n-1-2x)}{(n-1-2s)^2}\mbox{d}s.
\end{split}
\end{equation}
For $x\geq 0$ and $y<1/2$, we compute that
\begin{equation}\label{Usmall2}
\begin{split}
 &U(x,y,z)\\
&=U(\mbox{$\frac{1}{2}$,$\frac{1}{2}$},-y)\cdot \frac{2(n-1)(y-x)}{n-1+2y}+U(\mbox{$\frac{1}{2}$,$\frac{1}{2}$},z)\cdot \frac{2((n-1)x-z)}{n-1-2z}\\
&\quad +\int_{-1/2}^{-y} U(\mbox{$\frac{1}{2}$,$\frac{1}{2}$},s)\cdot \frac{2n(n-1)}{(n-1-2s)^2}\mbox{d}s+\int_{-y}^{z} U(\mbox{$\frac{1}{2}$,$\frac{1}{2}$},s)\cdot \frac{2(n-1)(1-2x)}{(n-1-2s)^2}\mbox{d}s.
\end{split}
\end{equation}
\end{cor}

The values of $U(\frac{1}{2},\frac{1}{2},s)$ can be extracted from Corollary \ref{cor+int}. 
In particular, the formula \eqref{Usmall2} can be applied for $x=y=z=0$, resulting in quite an involved, yet nonetheless explicit expression. 

\begin{rem}
It is straightforward to check that, for all $y>0$, the function $U$ satisfies the symmetry condition $U(-y,y,-y)=U(y,y,-y)$: compare the first and the fifth line in \eqref{form1}, see also \eqref{Usmall1} and \eqref{Usmall2}. This is in perfect consistence with the jump property of $X$ described at the end of Step 1. Indeed, as we noted there, when the left limit $X_{t-}$ is equal to $-Y_t$, then at time $t$ the process $X$ jumps from $-Y_t$ to $Y_t$. In other words, the points $(-y,y,-y)$ and $(y,y,-y)$ in the state space correspond to the same value of $U$.
\end{rem}

\section{Proof of Theorem \ref{main_theorem}} \label{sec4}
We now will prove that the function $U$ constructed in the previous section is indeed the value function of the optimal stopping problem \eqref{optimal_stopping} under consideration. We start with the majorization property.

\begin{lemma}
We have $U(x,y,z)\geq y-z-x^2$ for all $(x,y,z)$.
\end{lemma}
\begin{proof}
Suppose first that $y\geq \frac{1}{2}$ and $z\leq \varphi(y-\frac{1}{2})$. According to \eqref{form1}, we need to consider five cases. For $x\leq \varphi^{-1}(z)\leq 0$, the desired estimate is equivalent to $(2x-z-1)^2\geq 0$. If $x\leq 0\leq \varphi^{-1}(z)$, then
\begin{align*}
 U(x,y,z)&=y-z-\frac{2\varphi^{-1}(z)}{n-1}x+\big(\varphi^{-1}(z)\big)^2\geq y-z\geq y-z-x^2.
\end{align*}
For $0\leq x\leq \varphi^{-1}(z)$, the claim reads $(x-\varphi^{-1}(z))^2\geq 0$. If $\varphi^{-1}(z)\leq x\leq y-\frac{1}{2}$, then the majorization is actually an equality. Finally, for $x\geq y-\frac{1}{2}$, the desired bound becomes  $\left(x-y+\frac{1}{2}\right)^2\geq 0$, which is also trivial.

Now, suppose that $y<\frac{1}{2}$ or $z>\varphi(y-\frac{1}{2})$. It follows directly from \eqref{uy}, \eqref{Usmall1} and \eqref{Usmall2} that $U_y(x,y+,z)\leq 1$, that is, 
$$ \frac{\partial}{\partial y+}\big(U(x,y,z)-(y-z-x^2)\big)\leq 0.$$
Therefore, the majorization follows at once from the above analysis: we have
$$ U(x,y,z)-(y-z-x^2)\geq U(x,\varphi^{-1}(z)+\mbox{$\frac{1}{2}$},z)-(\varphi^{-1}(z)+\mbox{$\frac{1}{2}$}-z-x^2)\geq 0. \qedhere$$
\end{proof}

\begin{lemma}\label{Rem_symmetry}
For any $(x,y,z)$ and any bounded stopping time $\tau$ we have
\begin{equation}\label{main}
 \E_{x,y,z}(Y_\tau-Z_{\tau}-X_{\tau}^2)\leq U(x,y,z).
\end{equation}
\end{lemma}
\begin{proof} Roughly speaking, the argument rests on It\^o's formula and the majorization established in the previous section. However, since the function $U$ is not in $C^2$, there are some technical obstacles, which will be handled by an appropriate stopping procedure. We split the argument into a few parts. 

\smallskip

\emph{Part 1.} By continuity, we may and do assume that $y>0$ and $-y<z<0$. 
We introduce the increasing sequence $(\sigma_n)_{n\geq 0}$ of stopping times given inductively by
$$ \sigma_0=\inf\{t\geq 0:Y_t\geq \mbox{$\frac{1}{2}$}\}$$
 and, for $n\geq 0$, 
$$ \sigma_{2n+1}=\inf\{t\geq \sigma_{2n}:X_t=0\}, \qquad \sigma_{2n+2}=\inf\{t\geq \sigma_{2n+1}: X_t\in\{Y_t,Z_t\}\}.$$
It is straightforward to see that $\lim_{n\to \infty}\sigma_n=\infty$ almost surely. The function $U$ is of class $C^\infty$ on $$D\cap \{(x,y,z):x\neq 0\mbox{ and }y<1/2\}$$ and satisfies $ U_{xx}(x,y,z)=0$ on this set. Furthermore, we may easily check that $$(n-1)U_{x}(0-,y,z)=U_x(0+,y,z)$$ and  $$U_y(y,y,z)=U_z(z,y,z)=0$$ for all $y,\,z$, and that $$U(-y,y,-y)=U(y,y,-y)$$ for all $y$ (as noted in Remark \ref{Rem_symmetry} above). Therefore, by It\^o-Tanaka formula, we obtain
$$ \E_{x,y,z}U(X_{\tau\wedge \sigma_0},Y_{\tau\wedge \sigma_0},Z_{\tau\wedge \sigma_0})=U(x,y,z).$$
Note that the symmetry condition $U(-y,y,-y)=U(y,y,-y)$ guarantees that the jumps of $X$ do not contribute. 

\smallskip

\emph{Part 2.} If $X_{\tau\wedge \sigma_0}> 0$, then we look at the restriction $U^+=U|_{\mathcal{D}^+}$, where $\mathcal{D}^+=\mathcal{D}\cap \{x\leq 0\}$. Then $U^+(\cdot,y_0,z_0)$ is concave for any $y_0,\,z_0$ and satisfies $U^+_y(y_0,y_0,z_0)=0$. It\^o's formula then gives 
\begin{align*}
\E_{x,y,z}\big[U(X_{\tau\wedge \sigma_1},Y_{\tau\wedge \sigma_1},Z_{\tau\wedge \sigma_1})|\F_{\tau\wedge \sigma_0}\big]&=
\E_{x,y,z}\big[U^+(X_{\tau\wedge \sigma_1},Y_{\tau\wedge \sigma_1},Z_{\tau\wedge \sigma_1})|\F_{\tau\wedge \sigma_0}\big]\\
&\leq U^+(X_{\tau\wedge \sigma_0},Y_{\tau\wedge \sigma_0},Z_{\tau\wedge \sigma_0})\\
&= U(X_{\tau\wedge \sigma_0},Y_{\tau\wedge \sigma_0},Z_{\tau\wedge \sigma_0}).
\end{align*}

If $X_{\tau\wedge \sigma_0}< 0$ (which happens only if $y\geq 1/2$ and $x\leq 0$), we proceed similarly. Consider the restriction $U^-=U|_{\mathcal{D}^-}$, where $\mathcal{D}^-=\mathcal{D}\cap \{x\leq 0\}$. Then $U^-(\cdot,y_0,z_0)$ is concave for any $y_0,\,z_0$, satisfies $U^-_y(z_0,y_0,z_0)=0$ and $$U^-(-y_0,y_0,-y_0)=U^+(y_0,y_0,y_0).$$ Therefore, applying It\^o's formula (and noting that the latter identity allow to ignore the jumps of $X$), we obtain again that
\begin{align*}
\E_{x,y,z}\big[U(X_{\tau\wedge \sigma_1},Y_{\tau\wedge \sigma_1},Z_{\tau\wedge \sigma_1})|\F_{\tau\wedge \sigma_0}\big]
\leq  U(X_{\tau\wedge \sigma_0},Y_{\tau\wedge \sigma_0},Z_{\tau\wedge \sigma_0}).
\end{align*}
Consequently, we have shown that
$$ \E_{x,y,z}U(X_{\tau\wedge \sigma_1},Y_{\tau\wedge \sigma_1},Z_{\tau\wedge \sigma_1})=\E_{x,y,z}U(X_{\tau\wedge \sigma_0},Y_{\tau\wedge \sigma_0},Z_{\tau\wedge \sigma_0}).$$

\smallskip

\emph{Part 3.} 
Next, on the time interval $[\sigma_1,\sigma_2]$, the processes $Y$ and $Z$ remain unchanged. Therefore, since the function $U(\cdot,y_0,z_0)$ is concave and satisfies $$U_x(0-,y_0,z_0)=(n-1)U_x(0+,y_0,z_0),$$ applying the It\^o-Tanaka formula yields
$$ \E_{x,y,z}U(X_{\tau\wedge \sigma_2},Y_{\tau\wedge \sigma_2},Z_{\tau\wedge \sigma_2})\leq \E_{x,y,z}U(X_{\tau\wedge \sigma_1},Y_{\tau\wedge \sigma_1},Z_{\tau\wedge \sigma_1}).$$

\smallskip

\emph{Part 4.} Iterating the above argument, we see that the sequence  
$$\big(\E_{x,y,z}U(X_{\tau\wedge \sigma_n},Y_{\tau\wedge \sigma_n},Z_{\tau\wedge \sigma_n})\big)_{n\geq 0}$$ 
is nonincreasing and hence in particular $ \E_{x,y,z}U(X_{\tau\wedge \sigma_n},Y_{\tau\wedge \sigma_n},Z_{\tau\wedge \sigma_n})\leq U(x,y,z).$ 
By the previous lemma, this implies
$$ \E_{x,y,z}(Y_{\tau\wedge \sigma_n}-Z_{\tau\wedge \sigma_n})\leq U(x,y,z)+\E_{x,y,z}X_{\tau\wedge \sigma_n}^2\leq U(x,y,z)+\E_{x,y,z}X_{\tau}^2,$$
where in the last passage we have exploited the submartingale property of $X^2$. 
Letting $n\to \infty$, we see that the left-hand side tends to $\E_{x,y,z}(Y_{\tau}-Z_{\tau})$, by Lebesgue's monotone convergence theorem. This yields the desired claim.
\end{proof}

Together with the reasoning from the previous section, Lemma \ref{Rem_symmetry} identifies the explicit formula for the value function $\mathbb{U}$ of the optimal stopping process \eqref{stop_n}. Corollary \ref{cor1} thus follows.

\begin{cor} \label{cor1}
We have $\mathbb{U}=U$ on $D$.
\end{cor}

We are ready for the proof of our main result.

\begin{proof}[Proof of Theorem \ref{main_theorem}]
Let us write $\mathbb{P}$ instead of $\mathbb{P}_{0,0,0}$. By Lemma \ref{Rem_symmetry}, for any stopping time $\tau$ we have
$$ \E (Y_\tau-Z_\tau)\leq U(0,0,0)+\E X_\tau^2=U(0,0,0)+\E \tau.$$
Let us apply a scaling argument. As we have already discussed above, for any $\lambda>0$, $(\tilde{X}^\lambda_t)_{t\geq 0}=(\lambda X_{t/\lambda^2})_{t\geq 0}$ is a spider process and $\tilde{\tau}^\lambda=\tau \lambda^2$ is a stopping time of $\tilde{X}$. Consequently, applying the above inequality to $\tilde{X}^\lambda$, we obtain
$$ \E (Y_\tau-Z_\tau)=\lambda^{-1}\E (\tilde{Y}^\lambda_{\tilde{\tau}^\lambda}-\tilde{Z}^\lambda_{\tilde{\tau}^\lambda})\leq \lambda^{-1}U(0,0,0)+\lambda^{-1}\E \tilde{\tau}^\lambda=\lambda^{-1}U(0,0,0)+\lambda \E \tau.$$
Optimizing over $\lambda$, we get
$$ \E (Y_\tau-Z_\tau)\leq 2\sqrt{U(0,0,0)\E \tau},$$
which is the desired inequality. The equality is attained for the special stopping time $\tau$ considered in the previous section: indeed, as we have shown there, this particular $\tau$ satisfies $ \E(Y_\tau-Z_\tau)=U(0,0,0)+\E \tau$, and the right hand side is at least $ 2\sqrt{U(0,0,0)\E \tau}. $ This concludes the proof. 
\end{proof}

\end{document}